 \documentclass[9pt,journal]{IEEEtran} % ,draftcls,onecolumn
\newtheorem{theorem}{Theorem}

\newtheorem{lemma}[theorem]{Lemma}

\newtheorem{remark}{Remark}

\newcommand{\enma}[1]   {\ensuremath{#1}}

\newcommand{\non}{\nonumber}

\newcommand{\beq}{\begin{equation}}
\newcommand{\eeq}{\end{equation}}
\newcommand{\bseq}{\begin{subequations}}
\newcommand{\eseq}{\end{subequations}}
\newcommand{\beqn}{\begin{eqnarray}}
\newcommand{\eeqn}{\end{eqnarray}}
\newcommand{\ba}{\begin{array}}
\newcommand{\ea}{\end{array}}
\newcommand{\bct}{\begin{center}}
\newcommand{\ect}{\end{center}}
\newcommand{\btmz}{\begin{itemize}}
\newcommand{\etmz}{\end{itemize}}
\newcommand{\benum}{\begin{enumerate}}
\newcommand{\eenum}{\end{enumerate}}

\newcommand{\cH}{\enma{\mathcal H}}
\newcommand{\cL}{\enma{\mathcal L}}

\newcommand{\norm}[1]{\| #1 \|}                 %does not make large \|

\newcommand{\Ht}{ \mathcal{H}_{2}}

\newcommand{\inner}[2]{\left\langle #1,#2 \right\rangle}

\newcommand{\matbegin}{
        \left[
}
\newcommand{\matend}{
        \right]
}

\newcommand{\tbo}[2]{
  \matbegin \begin{array}{c}
       #1 \\ #2
       \end{array} \matend }

\newcommand{\tbt}[4]{
  \matbegin \begin{array}{cc}
       #1 & #2 \\ #3 & #4
       \end{array} \matend }

\newcommand{\be}{\begin{equation}}
\newcommand{\ee}{\end{equation}}

\newcommand{\cplxs}{ C\kern -.35em \rule{0.03 em}{.7 ex}~   }

\def\complex{\hbox{C\kern -.45em \rule{0.03 em}{1.5 ex}}~}

\newcommand{\bi}{\begin{itemize}}
\newcommand{\ei}{\end{itemize}}

\usepackage{amssymb,latexsym,color,amsmath,epsfig,graphicx,subfig,array}
\usepackage{setspace,cite,dsfont} 
\usepackage{cite}
\usepackage{algorithm,algorithmic}
\usepackage{hyperref}
\usepackage[usenames,dvipsnames,svgnames,table]{xcolor}
\usepackage{tikz}
\usetikzlibrary{calc,arrows,automata,positioning}
\usepackage{pgfplots}
\usepgfplotslibrary{patchplots}

\usepackage{listings}
\newcommand{\cC}{{\cal C}}
\newcommand{\cE}{{\cal E}}

\newcommand{\cP}{{\cal P}}

\newcommand{\cS}{{\cal S}}
\newcommand{\cT}{{\cal T}}

\newcommand{\mre}{\mathrm{e}}
\newcommand{\bbR}{\mathbb{R}}
\newcommand{\ds}{\displaystyle}

\newcommand{\tc}{\textcolor}

\newcommand{\tcr}{\textcolor{red}}

\DeclareMathOperator*{\argmin}{argmin}
\DeclareMathOperator*{\minimize}{minimize}
\DeclareMathOperator*{\subject}{subject~to}
\DeclareMathOperator*{\vect}{vec}

\newcommand{\mrj}{\mathrm{j}}

\newcommand{\ty}{{\tilde{y}}}

\newcommand{\DefinedAs}[0]{\mathrel{\mathop:}=}
\newcommand{\AsDefined}[0]{=\mathrel{\mathop:}}

\newcommand{\sign}{\mathrm{sign}}

\newcommand{\tx}{{\tilde{x}}}

\newcommand{\prox}{\mathbf{prox}}

\newcommand{\one}{\mathds{1}}

\newcommand{\fx}{f}
\newcommand{\gz}{g}

\newcommand{\iter}{l}
\newcommand{\tm}{\tilde m}

\newcommand{\zms}{z_\mu^\star}

\usepackage{accents}

\definecolor{newyellow}{rgb}{0.929,0.694,0.125}

\title{The proximal augmented Lagrangian method \\[-0.05cm] for nonsmooth composite optimization}
\author{Neil K.\ Dhingra, Sei Zhen Khong, and Mihailo R.\ Jovanovi\'c
\thanks{Financial support from under NSF Awards ECCS-1739210 and CNS-1544887 and AFOSR Award FA9550-16-1-0009 is gratefully acknowledged.}
\thanks{N.\ K.\ Dhingra is with Numerica Corporation, Fort Collins, CO 80528; S. Z. Khong is with the Department of Electrical and Electronic Engineering, University of Hong Kong, Pokfulam, Hong Kong, China; M.\ R.\ Jovanovi\'c is with the Ming Hsieh Department of Electrical Engineering, University of Southern California, Los Angeles, CA, 90089. E-mails: dhin0008@umn.edu, szkhong@hku.hk, mihailo@usc.edu.}}

\begin{document}
\maketitle
\vspace{-5ex}
	\begin{abstract}
We study a class of optimization problems in which the objective function is given by the sum of a differentiable but possibly nonconvex component and a nondifferentiable convex regularization term. We introduce an auxiliary variable to separate the objective function components and utilize the Moreau envelope of the regularization term to derive the proximal augmented Lagrangian -- a continuously differentiable function obtained by constraining the augmented Lagrangian to the manifold that corresponds to the explicit minimization over the variable in the nonsmooth term. The continuous differentiability of this function with respect to both primal and dual variables allows us to leverage the method of multipliers (MM) to compute optimal primal-dual pairs by solving a sequence of differentiable problems. The MM algorithm is applicable to a broader class of problems than proximal gradient methods and it has stronger convergence guarantees and a more refined step-size update rules than the alternating direction method of multipliers. These features make it an attractive option for solving structured optimal control problems. We also develop an algorithm based on the primal-descent dual-ascent gradient method and prove global (exponential) asymptotic stability when the differentiable component of the objective function is (strongly) convex and the regularization term is convex. Finally, we identify classes of problems for which the primal-dual gradient flow dynamics are convenient for distributed implementation and compare/contrast our framework to the existing approaches.
	\end{abstract}
	
	\vspace*{-3ex}
\section{Introduction}

	We study a class of composite optimization problems in which the objective function is a sum of a differentiable  but possibly nonconvex component and a convex nondifferentiable component. Problems of this form are encountered in diverse fields including compressive sensing~\cite{don06}, machine learning~\cite{hastibfri09}, statistics~\cite{biclev08}, image processing~\cite{golosh09}, and control~\cite{linfarjovTAC13admm}. In feedback synthesis, they typically arise when a traditional performance metric (such as the ${\cal H}_2$ or ${\cal H}_\infty$ norm) is augmented with a regularization function to promote certain structural properties in the optimal controller. For example, the $\ell_1$ norm and the nuclear norm are commonly used nonsmooth convex regularizers that encourage sparse and low-rank optimal solutions, respectively.

	The lack of a differentiable objective function precludes the use of standard descent methods for smooth optimization. Proximal gradient methods~\cite{parboy13} and their accelerated variants~\cite{bectebSIAM09} generalize gradient descent, but typically require the nonsmooth term to be separable over the optimization variable. Furthermore, standard acceleration techniques are not well-suited for problems with constraint sets that do not admit an easy projection (e.g., closed-loop stability).

	An alternative approach is to split the smooth and nonsmooth components in the objective function over separate variables which are coupled via an equality constraint. Such a reformulation facilitates the use of the alternating direction method of multipliers (ADMM)~\cite{boyparchupeleck11}. This augmented-Lagrangian-based method splits the optimization problem into subproblems which are either smooth or easy to solve. It also allows for a broader class of regularizers than proximal gradient and it is convenient for distributed implementation. However, there are limited convergence guarantees for nonconvex problems and parameter tuning greatly affects its convergence rate.
	
	The method of multipliers (MM) is the most widely used algorithm for solving constrained nonlinear programing problems~\cite{ber82,ber99,nocwri06}. In contrast to ADMM, it is guaranteed to converge for nonconvex problems and there are systematic ways to adjust algorithmic parameters. However, MM is not a splitting method and it requires {\em joint\/} minimization of the augmented Lagrangian with respect to {\em all\/} primal optimization variables. This subproblem is typically nonsmooth and as difficult to solve as the original optimization problem.

	To make this difficult subproblem tractable, we transform the augmented Lagrangian into the continuously differentiable proximal augmented Lagrangian by exploiting the structure of proximal operators associated with nonsmooth regularizers. This new form is obtained by constraining the augmented Lagrangian to the manifold that corresponds to the explicit minimization over the variable in the nonsmooth term. The resulting expression is given in terms of the Moreau envelope of the nonsmooth regularizer and is continuously differentiable. This allows us to take advantage of standard optimization tools, including gradient descent and quasi-Newton methods, and enjoy the convergence guarantees of standard~MM.
	
The proximal augmented Lagrangian also enables Arrow-Hurwicz-Uzawa primal-dual gradient flow dynamics. Such dynamics can be used to identify saddle points of the Lagrangian~\cite{arrhuruza59} and have enjoyed recent renewed interest in the context of networked optimization because, in many cases, the gradient can be computed in a distributed manner~\cite{waneli11}. Our approach yields a dynamical system with a continuous right-hand side for a broad class of nonsmooth optimization problems. This is in contrast to existing techniques which employ subgradient methods~\cite{nedozdTAC09} or use discontinuous projected dynamics~\cite{feipagAUT10,chemalcorSCL16,chemallowcor18} to handle inequality constraints. Furthermore, since the proximal augmented Lagrangian is not strictly convex-concave we make additional developments relative to~\cite{cheghacor17} to show asymptotic convergence. Finally, inspired by recent advances~\cite{lesrecpac16,husei16}, we employ the theory of integral quadratic constraints~\cite{megran97} to prove global exponential stability when the differentiable component of the objective function is strongly convex with a \mbox{Lipschitz continuous gradient.}

The rest of the paper is structured as follows. In Section~\ref{sec.problem}, we formulate the nonsmooth composite optimization problem {and} provide a brief background on proximal operators. In Section~\ref{sec.pal}, we exploit the structure of proximal operators to introduce the proximal augmented Lagrangian. In Section~\ref{sec.MM}, we provide an efficient algorithmic implementation of the method of multipliers using the proximal augmented Lagrangian. In Section~\ref{sec.pdda}, we prove global (exponential) asymptotic stability of primal-descent dual-ascent gradient flow dynamics under a (strong) convexity assumption. In Section~\ref{sec.ex}, we use the problems of edge addition in directed consensus networks {and optimal placement} to illustrate the effectiveness of our approach. We close the paper in Section~\ref{sec.remarks} with concluding remarks.

	\vspace*{-2ex}
\section{Problem formulation {and background}}
\label{sec.problem}

We consider a composite optimization problem,
\be
	\label{pr}
    \ba{rcl}
	\ds\minimize_x & \fx(x) \;+\; \gz\left(\cT(x)\right)
    \ea
\ee
where the optimization variable $x$ belongs to a finite-dimensional Hilbert space (e.g., $\bbR^n$ or $\bbR^{m \times n}$) equipped with an inner product $\inner{\cdot}{\cdot}$ and associated norm $\norm{\cdot}$. The function $\fx$ is continuously differentiable but possibly nonconvex, the function $g$ is convex but potentially nondifferentiable, and $\cT$ is a bounded linear operator. We  {further} assume that {$g$ is proper and lower semicontinuous, that}~\eqref{pr} is feasible, and that its minimum is finite.

Problem~\eqref{pr} is often encountered in structured controller design~\cite{jovdhiEJC16,wujovSCL17,mogjovTCNS17}, where $f$ is a measure of closed-loop performance, e.g., the $\cH_2$ norm, and the regularization term $g$ is introduced to promote certain structural properties of $\cT(x)$. For example, in wide-area control of power systems, $f$ measures the quality of synchronization between different generators and $g$ penalizes the amount of communication between them~\cite{dorjovchebulACC13,dorjovchebulTPS14,wudorjovTPS16}.

In particular, for $z \DefinedAs \cT (x) \in \bbR^m$, the $\ell_1$ norm, 
	$	
\norm{z}_1
	\DefinedAs 
		\sum |z_{i}|,
		$
is a commonly used convex proxy for promoting sparsity of $z$. For $z \in \bbR^{m \times  n}$, the nuclear norm,
$
    \norm{z}_{*}
    \DefinedAs
    \sum \sigma_i(z),
$
can be used to obtain low-rank solutions to~\eqref{pr}, where $\sigma_i(z)$ is the $i$th singular value. The indicator function,
$
    I_{\cC} (z)
    \DefinedAs
	\{
    0, \, z \in \cC;
    \,
    \infty, \, z \not\in \cC
	\}
$
associated with the convex set $\cC$ is the proper regularizer for enforcing $z \in \cC$.

Regularization of $\cT(x)$ instead of $x$ is important in the situations where the desired structure has a simple characterization in the co-domain of $\cT$. For example, such problems arise in spatially-invariant systems, where it is convenient to perform standard control design in the spatial frequency domain~\cite{bampagdah02} but necessary to promote structure in the physical space, and in consensus/synchronization networks, where the objective function is expressed in terms of the deviation of node values from the network average but it is desired to impose structure on the network edge weights~\cite{wujovSCL17,mogjovTCNS17}.

	\vspace*{-2ex}
\subsection{{Background on proximal operators}}
\label{sec.proxop}

Problem~\eqref{pr} is difficult to solve directly because $f$ is, in general, a nonconvex function and $g$ is typically not differentiable. Since the existing approaches and our method utilize proximal operators, we first provide a brief overview; for additional information, see~\cite{parboy13}.

The proximal operator of the function $g$ is given by
\begin{subequations}
	\beq
	\label{eq.prox}
	\prox_{\mu g} (v)
	\; \DefinedAs \;
	\argmin\limits_{x}
	\,
	\left(
	g (x)	
	\; + \;
	\tfrac{1}{2\mu}
	\,
	\norm{x \, - \, v}^2
	\right)
	\eeq
and the associated optimal value specifies its Moreau envelope,	
	\beq
	M_{\mu g} (v)
	\; \DefinedAs \;
	g (\prox_{\mu g} (v))	
	\; + \;
	\tfrac{1}{2\mu}
	\,
	\norm{\prox_{\mu g} (v) \, - \, v}^2
	\eeq
where $\mu > 0$. The Moreau envelope is a continuously differentiable function, even when $g$ is not, and its gradient{~\cite{parboy13}} is given by
\be
	\nabla M_{\mu g}(v)
	\;=\;
	\tfrac{1}{\mu}
	\!
	\left(v \;-\; \prox_{\mu g}(v) \right).
	\label{eq.gradM}
\ee
\end{subequations}
For example, when $g$ is the $\ell_1$ norm, 
	$g (z) = \norm{z}_1 =\sum |z_i|$,
the proximal operator is determined by soft-thresholding,
	$
    	\prox_{\mu g}(v_i)
	=
	\cS_{\mu} (v_i)
	\DefinedAs 
         \sign (v_i)
         \max 
         \, 
         ( |v_i| - \mu, \, 0 ),
         $
the associated Moreau envelope is the Huber function,
	$
    M_{\mu g}(v_i)
    =
    \{
    \tfrac{1}{2\mu} \, v_i^2,
    \, 
    |v_i|
    \leq
    \mu;
    \, 
    |v_i| - \tfrac{\mu}{2},
    \,
    |v_i|
    \geq
    \mu
    \},
    $
and the gradient of this Moreau envelope is the saturation function,
	$
    \nabla M_{\mu g}(v_i)
    =
    \sign(v_i) \min \, ( |v_i|/\mu, \,1 ).
	$

	\vspace*{-2ex}
\subsection{Existing algorithms}

\subsubsection{Proximal gradient}

The proximal gradient method generalizes standard gradient descent to certain classes of nonsmooth optimization problems. This method can be used to solve~\eqref{pr} when $g(\cT)$ has an easily computable proximal operator. When $\cT = I$, the proximal gradient method for problem~\eqref{pr} with step-size $\alpha_l$ is given by,
	\[
	\ba{r}
	x^{\iter+1}
	~ = ~
	\prox_{\alpha_\iter g}
	( x^\iter \,-\, \alpha_\iter\nabla f(x^\iter) ).
	\ea
	\]
When $\gz = 0$, the proximal gradient method simplifies to standard gradient descent, and when $\gz$ is indicator function of a convex set, it simplifies to projected gradient descent. The proximal gradient algorithm applied to the $\ell_1$-regularized least-squares problem (LASSO)
\be
    \minimize\limits_x
    ~\,
    \tfrac{1}{2} \, \norm{Ax \,-\, b}^2
    \;+\;
    \gamma
    \,
    \norm{x}_1
    \label{eq.lasso}
\ee
where $\gamma$ is a positive regularization parameter, yields the Iterative Soft-Thresholding Algorithm (ISTA)~\cite{bectebSIAM09},
	$
    x^{\iter+1}
    =
    \cS_{\gamma \alpha_\iter}
    (
    x^\iter
    -
    \alpha_\iter
    A^T(Ax^\iter - b)
    ).
	$
This method is effective only when the proximal operator of $g(\cT)$ is easy to evaluate.
Except in special cases, e.g, when $\cT$ {is diagonal}, efficient computation of $\prox_{\mu g(\cT)}$ does not necessarily follow from an efficiently computable $\prox_{\mu g}$. This makes the use of proximal gradient method challenging for many applications and its convergence can be slow. Acceleration techniques improve the convergence rate~\cite{bectebSIAM09,lilin15}, but they do not perform well in the face of constraints such as closed-loop stability.

\subsubsection{Augmented Lagrangian methods}

A common approach for dealing with a nondiagonal linear operator $\cT$ in~\eqref{pr} is to introduce an additional optimization variable $z$ 
    \be
    \label{pr.split}
    \ba{ll}
    \minimize\limits_{x, \, z}
    & f(x)
    \; + \;
    g(z)
    \\[0.15cm]
    \text{subject to}
    &
    \cT(x)
    \; - \;
    z
    \; = \;
    0.
    \ea
    \ee
The augmented Lagrangian is obtained by adding a quadratic penalty on the violation of the linear constraint to the regular Lagrangian associated with~\eqref{pr.split},
    \[
    \cL_\mu(x,z;y)
	\, = \,
    f(x)
    \, + \,
    g(z)
    \, + \,
    \inner{y}{\cT(x) \, - \, z}
     \,+\,
    \tfrac{1}{2\mu}
    \,
    \| \cT(x) \,-\, z \|^2
    \]
where $y$ is the Lagrange multiplier and $\mu$ is a positive parameter.

ADMM solves~\eqref{pr.split} via an iteration which involves minimization of $\cL_\mu{(x,z;y)}$ separately over $x$ and $z$ and a gradient ascent update (with step-size $1/\mu$) of $y$~\cite{boyparchupeleck11},
\begin{subequations}
	\label{eq.ADMM}
	\begin{IEEEeqnarray}{rcl}
	x^{k+1}
	&\;~=~\;&
	\argmin\limits_x
	\,
	\cL_\mu(x, \, z^k; \, y^k)
	\label{eq.ADMMx}
	\\[-0.1cm]
	z^{k+1}
	&\;~=~\;&
	 \argmin\limits_z
	 \,
	 \cL_\mu(x^{k+1}, \, z; \, y^k)
	 \label{eq.ADMMz}
	\\[-0.1cm]
	y^{k+1}
	&\;~=~\;&
	y^k
	~+~
	\tfrac{1}{\mu} \, (\cT(x^{k+1}) \, - \, z^{k+1}).
	\label{eq.ADMMy}
\end{IEEEeqnarray}
\end{subequations}
{ADMM is appealing because, even when $\cT$ is nondiagonal,} the $z$-minimization step amounts to {evaluating $\prox_{\mu g}$, and} the $x$-minimization step amounts to solving a smooth (but possibly nonconvex) optimization problem. Although it was recently shown that ADMM is guaranteed to converge {to a stationary point of~\eqref{pr.split} for some classes of nonconvex problems~\cite{honluorazSIAM16}}, its rate of convergence is strongly influenced by the choice of $\mu$.

The method of multipliers (MM) is the most widely used algorithm for solving constrained nonconvex optimization problems~\cite{ber82,congoutoi91} and it guarantees convergence to a local minimum. In contrast to ADMM, each MM iteration requires {\em joint\/} minimization of the augmented Lagrangian with respect to the primal variables $x$ and $z$,
\begin{subequations}
	\begin{IEEEeqnarray}{rcl}
	(x^{k+1},z^{k+1})
	&\;~=~\;&
	\argmin\limits_{x, \, z}
	\,
	\cL_\mu(x, \, z; \, y^k)
	\label{eq.MMxz}
    \\[-0.1cm]
    y^{k+1}
	&\;~=~\;&
    y^k
    ~+~
    \tfrac{1}{\mu} \, (\cT(x^{k+1}) \,-\, z^{k+1}).
    	\label{eq.MMy}
	\end{IEEEeqnarray}
	\label{eq.MM}
	\end{subequations}
{It is possible to refine MM to allow for inexact solutions to the $(x,z)$-minimization subproblem and adaptive updates of the penalty parameter $\mu$. However, until now, MM has not been a feasible choice for solving~\eqref{pr.split} because the nonconvex and nondifferentiable $(x,z)$-minimization subproblem is as difficult as the original problem~\eqref{pr}.}

	\vspace*{-1ex}
\section{The proximal augmented Lagrangian}
	\label{sec.pal}

{We next derive the proximal augmented Lagrangian, a continuously differentiable function resulting from explicit minimization of the augmented Lagrangian over the auxiliary variable $z$.} This brings the $(x,z)$-minimization problem~\eqref{eq.MMxz} into a form that is continuously differentiable with respect to both $x$ and $y$ and facilitates the use of a wide suite of standard optimization tools for solving~\eqref{pr}. In particular, as described below, our approach enables the method of multipliers and the Arrow-Hurwicz-Uzawa gradient \mbox{flow dynamics method.} 

	\vspace*{-2ex}
\subsection{Derivation of the proximal augmented Lagrangian}

The first main result of the paper is provided in Theorem~\ref{thm.diff}. We use the proximal operator of the function $g$ to eliminate the {auxiliary} optimization variable $z$ from the augmented Lagrangian and transform~\eqref{eq.MMxz} into a tractable continuously differentiable problem.
\begin{theorem} \label{thm.diff}
For a proper, lower semicontinuous, and convex function function $g$, minimization of the augmented Lagrangian $\cL_\mu(x,z; y)$ associated with problem~\eqref{pr.split} over $(x,z)$ is equivalent to minimization of the {\em proximal augmented Lagrangian\/}
\be
	\cL_\mu(x;y)
	\;\DefinedAs\;
	f(x)
    \; + \;
    M_{\mu g}
    (\cT(x) \, + \, \mu y)
    \; - \;
    \tfrac{\mu}{2} \, \| y \|^2
    \label{eq.alprox}
\ee
over $x$. Moreover, if $f$ is continuously differentiable $\cL_\mu(x;y)$ is continuously differentiable over $x$ and $y$, and if $f$ has a Lipschitz continuous gradient $\nabla \cL_\mu(x;y)$ is Lipschitz continuous.
\end{theorem}
\begin{IEEEproof}
 Through the completion of squares, {the augmented Lagrangian $\cL_\mu$ associated with~\eqref{pr.split} can be equivalently written} as
    \[
    \cL_\mu(x,z;y)
    \, = \,
    f(x)
    \, + \,
    g(z)
    \, + \,
    \tfrac{1}{2\mu} \, \norm{z \, - \, (\cT(x) + \mu y)}^2
    \, - \,
    \tfrac{\mu}{2} \, \norm{y}^2.
    \]
Minimization with respect to $z$ yields an explicit expression,
	\be
	\zms(x,y)
	\; = \;
	\prox_{\mu g}(\cT(x) \, + \, \mu y)
	\label{eq.zstar}
	\ee
and substitution of $\zms$ into the augmented Lagrangian provides~\eqref{eq.alprox}, i.e., 
$
	\cL_\mu (x; y)
	=
	\cL_\mu(x, \zms(x,y); y).
$
Continuous differentiability of $\cL_\mu(x;y)$ follows from continuous differentiability of $M_{\mu g}$ and Lipschitz continuity of $\nabla \cL_\mu(x;y)$ follows from Lipschitz continuity of $\prox_{\mu g}$ and boundedness of the linear operator $\cT$; see~\eqref{eq.gradM}.
\end{IEEEproof}

Expression{~\eqref{eq.alprox}}, that we refer to as the {\em proximal augmented Lagrangian\/}, characterizes $\cL_\mu{(x,z;y)}$ on the manifold corresponding to explicit minimization over the {auxilary} variable $z$. {Theorem~\ref{thm.diff} allows} {\em joint\/} minimization of the augmented Lagrangian with respect to $x$ and $z$, which is in general a nondifferentiable problem, {to} be achieved by minimizing differentiable function~\eqref{eq.alprox} over $x$. It thus facilitates the use of the method of multipliers {in Section~\ref{sec.MM} and the Arrow-Hurwicz-Uzawa gradient flow dynamics in Section~\ref{sec.pdda}.}
	
	\begin{remark}
The proximal augmented Lagrangian can be derived even in the presence of a more general linear constraint,
	\begin{subequations}
	\be
	\ba{rl}
	\ds\minimize\limits_{x_1, \, x_2}
	&
	f(x_1)
	\;+\;
	g(x_2)
	\\[0.15cm]
	\subject
	&
	\cT_1 (x_1)
	\;+\;
	\cT_2 (x_2)
	\;=\;
	0.
	\ea
	\ee
Introduction of an additional auxiliary variable $z$ in the nonsmooth part of the objective function $g$, can be used to bring this two-block optimization problem into the following form,
	\be
	\ba{rl}
	\ds\minimize\limits_{x_1, \, x_2, \, z}
	&
	f(x_1)
	\;+\;
	g(z)
	\\[0.15cm]
	\subject
	&
	\cT_1 (x_1)
	\;+\;
	\cT_2 (x_2)	
	\;=\;
	0,
	~~
	x_2
	\; - \;
	z
	\;=\;
	0.
	\ea
	\label{eq.2block}
\ee
	\label{eq.more-general}
	\end{subequations}
\!\!Via an analogous procedure to that described in Theorem~\ref{thm.diff}, explicit minimization with respect to $z$ can be employed to eliminate it from the augmented Lagrangian and obtain a continuously differentiable function of both primal ($x_1,x_2$) and dual ($y_1,y_2$) variables,
\[
	\ba{lcr}
	{\cal L}_{\mu}(x_1,x_2;y_1,y_2)
	& \!\!\!\! = \!\!\!\! &
	f(x_1)
	\, + \,
	\tfrac{1}{2 \mu} \, \norm{\cT_1 (x_1)
	+
	\cT_2 (x_2) + \mu y_1}^2
	~ +
	\\[0.15cm]
	& \!\!\!\! \!\!\!\! &
	M_{\mu g}(x_2 + \mu y_2)
	\, - \,
	\tfrac{\mu}{2} \, \norm{y_1}^2
	\, - \,
	\tfrac{\mu}{2} \, \norm{y_2}^2.
	\ea
\]
Here, $y_1$ and $y_2$ are the Lagrange multipliers associated with the respective linear constraints in~\eqref{eq.2block} and, for simplicity, we use single parameter $\mu$ in the augmented Lagrangian. This approach has numerous advantages over standard ADMM; e.g., it can be readily extended to multi-block optimization problems for which ADMM is not guaranteed to converge in general~\cite{cheheyeyua16}. These extensions are outside of the scope of the present study and will be reported~elsewhere.	
	\end{remark}

	\vspace*{-2ex}
\subsection{MM using the proximal augmented Lagrangian} \label{sec.MM}

Theorem~\ref{thm.diff} allows us to solve nondifferentiable subproblem~\eqref{eq.MMxz} by minimizing the continuously differentiable proximal augmented Lagrangian $\cL_\mu(x; \, y^k)$ over $x$. We note that similar approach was also applied to MM in~\cite{roc76prox} for the particular case in which $g$ is the indicator function of a convex set. Relative to ADMM, our customized MM algorithm guarantees convergence to a local minimum and offers systematic update rules for the parameter $\mu$. Relative to proximal gradient, we can solve~\eqref{pr} with a general bounded linear operator $\cT$ and can incorporate second order information about $f$. 

Using reformulated expression~\eqref{eq.alprox} for the augmented Lagrangian, {MM} minimizes $\cL_{\mu} (x; y^k)$ over the primal variable $x$ and updates the dual variable $y$ using gradient ascent with step-size $1/\mu$,
	\begin{align}
	\tag*{(MMa)}
	\label{eq.MMa1}
	x^{k+1}
	\,&=\;
	\argmin\limits_{x}
	\,
	\cL_{\mu} (x; y^k)
	\\[-0.1cm]	
	y^{k+1}
	\,&=\;
	y^{k}
	\; + \;
		\tfrac{1}{\mu}
	 \,
	\nabla_y \; \! \cL_{\mu} (x^{k+1}; y^k)
	\tag*{(MMb)}
	\label{eq.MMa2}
	\end{align}
where
	$
	\nabla_y \; \! \cL_{\mu} (x^{k+1}; y^k)
	\DefinedAs 
	\cT (x^{k+1})
	- 
	z_\mu^{\star}(x^{k+1},y^k)
	= 
	\cT (x^{k+1})
	 - 
	\prox_{\mu g}(\cT(x^{k+1}) + \mu y^k)
	$
denotes the primal residual.

In contrast to ADMM, our approach does not attempt to avoid the lack of differentiability of $g$ by fixing $z$ to minimize over $x$. By constraining $\cL_\mu(x,z;y)$ to the manifold resulting from explicit minimization \mbox{over $z$}, we guarantee continuous differentiability of the proximal augmented Lagrangian $\cL_\mu(x;y)$. {MM is a gradient ascent algorithm on the Lagrange dual of a version of~\eqref{pr.split}, with the same constraint, in which the objective function is replaced by $f(x) + g(z) + \tfrac{1}{2\mu}\norm{\cT(x)-z}^2$; see~\cite[Section 2.3]{boyparchupeleck11} and~\cite{roc76}. Since its closed-form expression is typically unavailable, MM uses the $(x,z)$-minimization subproblem~\eqref{eq.MMxz} to evaluate this dual computationally and then take a gradient ascent step~\eqref{eq.MMy} in $y$. ADMM avoids this issue by solving simpler, separate subproblems over $x$ and $z$. However, the $x$ and $z$ minimization steps~\eqref{eq.ADMMx} and~\eqref{eq.ADMMz} do not solve~\eqref{eq.MMxz} and thus unlike the $y$-update~\eqref{eq.MMy} in MM, the $y$-update~\eqref{eq.ADMMy} in ADMM is not a gradient ascent step on the ``strengthened dual''. MM thus offers stronger convergence results~\cite{ber82,boyparchupeleck11} and may lead to fewer $y$-update steps.

	\begin{remark}
	The proximal augmented Lagrangian enables MM because the $x$-minimization subproblem in MM~\ref{eq.MMa1} is not more difficult than in ADMM~\eqref{eq.ADMMx}. For LASSO problem~\eqref{eq.lasso}, the $z$-update in ADMM~\eqref{eq.ADMMz} is given by soft-thresholding, $z^{k+1} = \cS_{\gamma \mu}(x^{k+1} + \mu y^k)$, and the $x$-update~\eqref{eq.ADMMx} requires minimization of the quadratic function~\cite{boyparchupeleck11}. In contrast, the $x$-update~\ref{eq.MMa1} in MM requires minimization of
	$
    	(1/2) \, \norm{Ax \, - \, b}^2
    	+
	M_{\mu_k g}
	(x \, + \, \mu_k y^k),
	$
where $M_{\mu_k g} (v)$ is the Moreau envelope associated {with the $\ell_1$ norm;} i.e., the Huber function. Although in this case the solution to~\eqref{eq.ADMMx} can be characterized explicitly by a matrix inversion, this is not true in general. The computational cost associated with solving either~\eqref{eq.ADMMx} or~\ref{eq.MMa1} using first-order methods scales at the same rate.
	\end{remark}
	
\subsubsection{Algorithm}

The procedure outlined in~\cite[Algorithm~17.4]{nocwri06} allows minimization subproblem~\ref{eq.MMa1} to be inexact, provides a method for adaptively adjusting $\mu_k$, and describes a more refined update of the Lagrange multiplier $y$. We incorporate these refinements into our proximal augmented Lagrangian algorithm for solving~\eqref{pr.split}. In Algorithm~\ref{alg.mm}, $\eta$ and $\omega$ are convergence tolerances, and $\mu_{\min}$ is a minimum value of the parameter $\mu$. {Because of the equivalence established in Theorem~\ref{thm.diff}}, convergence to a local minimum follows from the convergence {results for} the standard method of multipliers~\cite{nocwri06}.

\begin{algorithm}
\caption{MM using the proximal augmented Lagrangian.}
\label{alg.mm}
\begin{algorithmic}
\STATE \textbf{input:} Initial point $x^0$ and Lagrange multiplier $y^0$
\STATE \textbf{initialize:} $\mu_0 = 10^{-1}$, $\mu_{\min} = 10^{-5}$, $\omega_0 = \mu_0$, and $\eta_0 = \mu_0^{0.1}$
\vspace*{0.15cm}
\STATE \textbf{for} $k = 0, 1, 2, \ldots$
\\[.1cm]
~\,\quad Solve~\ref{eq.MMa1} such that~~~
	$
	\norm{\nabla_x \; \! \cL_\mu(x^{k+1},y^k)}
	\, \leq \,
	\omega_k
	$
\\[.1cm]
~\,\quad \textbf{if} $\norm{\nabla_y \; \! \cL_{\mu_k} (x^{k+1}; y^k)} \, \leq \, \eta_k$ 
\\[.1cm]
~\,\quad \quad \textbf{if} $\norm{\nabla_y \cL_{\mu_k} (x^{k+1}; y^k)} \leq \eta$  and \mbox{$\norm{\nabla_x  \cL_\mu(x^{k+1},y^k)}  \leq \omega$} 
\\[.1cm]
~\,\quad~\,\quad~\,\quad \textbf{stop} with approximate solution $x^{k+1}$
\\[.1cm]
~\,\quad \quad \textbf{else:}
\\[.1cm]
~\,\quad \quad \quad $y^{k+1} \, = \, y^k + \frac{1}{\mu_k} \nabla_y \; \! \cL_{\mu_k} (x^{k+1}; y^k)$, \hfill $\mu_{k+1} \, = \, \mu_k$ ~~~~~
\\[.1cm]
~\,\quad \quad \quad $\eta_{k+1} \, = \, \eta_k \, \mu_{k+1}^{0.9}$,  \hfill $\omega_{k+1} \, = \, \omega_k \, \mu_{k+1}$ ~~~~~
\\[.1cm]
~\,\quad \textbf{else:}
\\[.1cm]
~\,\quad \quad $y^{k+1} \, = \, y^k$, \hfill $\mu_{k+1} \, = \, \max\{\mu_k/5,\mu_{\min}\}$ ~~~
\\
~\,\quad  \quad $\eta_{k+1} \, = \,  \mu_{k+1}^{0.1}$, \hfill $\omega_{k+1} \, = \,  \mu_{k+1}$~~~~~
\end{algorithmic}
\end{algorithm}

\subsubsection{Minimization of $\cL_\mu (x;y)$ over $x$}
	
{MM alternates between minimization of $\cL_\mu (x;y)$ with respect to $x$ (for fixed values of $\mu$ and $y$) and the update of $\mu$ and $y$. Since $\cL_\mu(x; y)$ is once continuously differentiable, many techniques can be used to find a solution to subproblem~\ref{eq.MMa1}. We next summarize three of them.

\subsubsection*{Gradient descent}
The gradient with respect to $x$ is given by,
\[
	\nabla_x \cL_\mu(x; y)
	=
	\nabla f(x)
	\, + \,
	\tfrac{1}{\mu}\,
	\cT^\dagger(\cT(x) + \mu y \,-\, \prox_{\mu g}(\cT(x) + \mu y))
\]
where $\cT^\dagger$ is the adjoint of $\cT$, $\inner{z}{\cT(x)} = \inner{\cT^\dagger(z)}{x}$. 
{Backtracking conditions such as the Armijo rule can be used to select a step-size.}

\subsubsection*{Proximal gradient}
Gradient descent does not exploit the structure of the Moreau envelope of the function $g$; in some cases, it may be advantageous to use proximal operator associated with the Moreau envelope to solve~\ref{eq.MMa1}. In particular, when $\cT = I$,~\eqref{eq.prox} and~\eqref{eq.gradM} imply that $\prox_{\alpha M_{\mu g}}(v) = x^*$ where $x^*$ satisfies,
	$
	x^*
	=
	\tfrac{1}{\mu \, + \, \alpha}
	\left(
	\alpha
	\,
	\prox_{\mu g}(x^*)
	+
	\mu
	\,
	 v
	\right).
	$
If $g$ is separable and has an easily computable proximal operator, its Moreau envelope also has an easily computable proximal operator. In~\cite{dhijovACC16}, proximal gradient methods were used for subproblem~\ref{eq.MMa1} to solve a sparse feedback synthesis problem introduced in~\cite{linfarjovTAC13admm}. Computational savings were shown relative to standard proximal gradient method and ADMM.

\subsubsection*{Quasi-Newton method} \label{sec.qn}

Although $\prox_{\mu g}$ is typically not differentiable, it is Lipschitz continuous and therefore differentiable almost everywhere~\cite{rocwet09}. To improve computational efficiency, we employ the limited-memory Broyden-Fletcher-Goldfarb-Shanno (L-BFGS) method~\cite[Algorithm 7.4]{nocwri06} which estimates the Hessian ${\nabla_{xx} \cL_\mu(x;y^k)}$ using first-order information and is guaranteed to converge for convex functions with Lipschitz continuous gradients~\cite{bongillemsag13}.

\begin{remark}
For regularization functions that do not admit simply computable proximal operators, $\prox_{\mu g}$ has to be evaluated numerically by solving~\eqref{eq.prox}. If this is expensive, the primal-descent dual-ascent algorithm of Section~\ref{sec.pdda} offers an appealing alternative because it requires one evaluation of $\prox_{\mu g}$ per iteration. When the regularization function $g$ is nonconvex, the proximal operator may not be single-valued and the Moreau envelope may not be continuously differentiable. In spite of this, the convergence of proximal algorithms has been established for nonconvex, proper, lower semicontinuous regularizers that obey the Kurdyka-\L{}ojasiewicz inequality~\cite{bolsabteb14}.
\end{remark}

	\vspace*{-2ex}
\section{Arrow-Hurwicz-Uzawa gradient flow}
    \label{sec.pdda}
    	% \vspace*{-1ex}

We now consider an alternative approach to solving~\eqref{pr}. Instead of minimizing over the primal variable and performing gradient ascent in the dual, we simultaneously update the primal and dual variables to find the saddle point of the augmented Lagrangian. The continuous differentiability of $\cL_\mu (x; y)$ established in Theorem~\ref{thm.diff} enables the use of Arrow-Hurwicz-Uzawa gradient flow dynamics~\cite{arrhuruza59},
	 \be
	\tbo{\dot{x}}{\dot{y}}
	\;=\;
	\tbo{-\nabla_x \; \! \cL_\mu(x;y)}{\phantom{-}\nabla_y \; \! \cL_\mu(x;y)}.
	\label{eq.ctdyn}
	\tag*{(GF)}
\ee
In Section~\ref{sec.as}, we show that the gradient flow dynamics~\ref{eq.ctdyn} globally converge to the set of saddle points of the proximal augmented Lagrangian $\cL_\mu (x; y)$ for a convex $f$ with a Lipschitz continuous gradient. In Section~\ref{sec.IQCs}, we employ the theory of IQCs to establish global exponential stability for a strongly convex $f$ with a Lipschitz continuous gradient and estimate convergence rates. Finally, in Section~\ref{sec.di} we identify classes of problems for which dynamics~\ref{eq.ctdyn} are convenient for distributed implementation and compare/contrast \mbox{our framework to the existing approaches.}

	\vspace*{-3ex}
\subsection{Global asymptotic stability for convex $f$}
	\label{sec.as}

We first characterize the optimal primal-dual pairs of optimization problem~\eqref{pr.split} with the Lagrangian,
$
	f(x)
	+
	g(z)
	+
	\inner{y}{\cT(x) - z}.
$
The associated first-order optimality conditions are given by,
\begin{subequations}
    \label{pr.cond}
\begin{IEEEeqnarray}{rcl}
	0
	&\;~=~\;&
	\nabla f(x^\star)
	\;+\;
	\cT^\dagger(y^\star)
	\label{pr.condgrad}
	\\[-0.05cm]
	0
	&\;~\in~\;&
	\partial g(z^\star)
	\;-\;
	y^\star
	\label{pr.condsubgrad}
	\\[-0.05cm]
	0
	&\;~=~\;&
	\cT(x^\star)
	\;-\;
	z^\star
	\label{pr.condfeas}
\end{IEEEeqnarray}
\end{subequations}
where $\partial g$ is the subgradient of $g$. Clearly, these are equivalent to the optimality condition for~\eqref{pr}, {i.e.,} $0 \in \nabla f(x^\star) + \cT^\dagger(\partial g(\cT(x^\star)))$. 
Even though we state the result for $x \in \bbR^n$ and $\cT(x) = Tx$ where $T \in \bbR^{m \times n}$ is a given matrix, the proof for $x$ in a Hilbert space and a {bounded} linear operator $\cT$ follows from similar arguments.
		
\begin{theorem}
Let $f$ be a continuously differentiable convex function with a Lipschitz continuous gradient and let $g$ be a {proper, lower semicontinuous,} convex function. Then, the set of optimal primal-dual pairs $(x^\star,y^\star)$ of~\eqref{pr.split} for the gradient flow dynamics~\ref{eq.ctdyn},
\be
	\tbo{\dot{{x}}}{\dot{{y}}}
	\;=\;
	\tbo{-\left(\nabla f(x) \;+\; T^T\nabla M_{\mu g}(Tx \,+\, \mu y)\right)}
    {\mu \, ( \nabla M_{\mu g}(Tx \,+\, \mu y) \;-\; y )}
	\label{eq.ctdyn1}
	\tag*{(GF1)}
\ee
is globally asymptotically stable (GAS) and $x^\star$ is a solution of~\eqref{pr}. 
\end{theorem}

\begin{IEEEproof}
We introduce a change of variables $\tx \DefinedAs x - x^\star$, $\ty \DefinedAs y - y^\star$ and a Lyapunov function candidate,
	$
	V(\tx,\ty)
	=
    \tfrac{1}{2}\inner{\tx}{\tx}
	+
    \tfrac{1}{2}\inner{\ty}{\ty},
	$
where $(x^\star,z^\star;y^\star)$ is an optimal solution to~\eqref{pr.split} that satisfies~\eqref{pr.cond}. The dynamics in the $(\tilde x, \tilde y)$-coordinates are given by,
\be
	\tbo{\dot{\tilde{x}}}{\dot{\tilde{y}}}
	\;=\;
	\tbo{- ( \nabla f(x) \, - \, \nabla f(x^\star) \,+\, (1/\mu)\, T^T\tm )}
    	{\tm \,-\, \mu \, \tilde y}
\label{eq.tilde}
\ee
where $\tm = \mu \, 
	( 
	\nabla M_{\mu g}(Tx +  \mu y)
	-   
	\nabla M_{\mu g}(Tx^\star + \mu y^\star) 
	)$
can be expressed as
	\be
	\ba{rrl}
	\tm
	& \!\!\! \DefinedAs \!\!\! &
	\tilde{v} \, - \, \tilde{z}
	\\[0.1cm]
	\tilde{v}
	& \!\!\! \DefinedAs \!\!\! &
	T\tx \, + \, \mu \ty
	\, = \,
	(Tx + \mu y)
	\, - \,
	(Tx^\star + \mu y^\star)
	\\[0.1cm]
	\tilde{z}
	& \!\!\! \DefinedAs \!\!\! &
	\prox_{\mu g}(Tx + \mu y)
	\,-\,
	\prox_{\mu g}(Tx^\star + \mu y^\star).
	\ea
	\label{eq.mtilde}
	\ee
The derivative of $V$ along the solutions of~\eqref{eq.tilde} is given by
	\[
	\ba{rcl}
	\!
	\dot{V}
	& \!\!\!\! = \!\!\!\! &
	-\inner{\tx}{\nabla f(x) - \nabla f(x^\star)}
	- 
	\tfrac{1}{\mu} \, \norm{T \tx}^2
	+ 
	\tfrac{1}{\mu} \inner{T \tx - \mu \ty}{\tilde{z}}
	\\[0.15cm]
	\!
	& \!\!\!\! = \!\!\!\! &
	-\inner{\tx}{\nabla f(x) - \nabla f(x^\star)}
	- 
	\tfrac{1}{\mu} 
	\left(
	\norm{T \tx}^2
	- 
	2 \inner{T \tx}{\tilde{z}}
	+
	\inner{\tilde{v}}{\tilde{z}}
	\right).
	\ea
\]
Since $f$ is convex with an $L_f$-Lipschitz continuous gradient and since $\prox_{\mu g}$ is firmly nonexpansive~\cite{parboy13}, i.e., $ \inner{\tilde{v}}{\tilde{z}} \geq \norm{\tilde{z}}^2$, we have
	\be
	\dot{V} (\tx,\ty)
	\; \leq \; 
	- \tfrac{1}{L_f} \, \norm{\nabla f(x) \, - \, \nabla f(x^\star)}^2
	\; - \;
	\tfrac{1}{\mu} \, \norm{T \tx \, - \, \tilde{z}}^2.
	\label{eq.Vdot}
	\ee
Thus, $\dot{V} \leq 0$ and each point in the set of optimal primal-dual pairs $(x^\star,y^\star)$ is stable in the sense of Lyapunov. 

The right-hand-side in~\eqref{eq.Vdot} becomes zero when $\nabla f(x) = \nabla f(x^\star)$ and $T \tx = \tilde{z}$. Under these conditions, we have $\dot{V} = -\inner{T^T \ty}{\tx}$ and the set of points for which $\dot{V} = 0$ is given by ${\cal D} = \{ (x,y); \nabla f(x) = \nabla f(x^\star)$, $T \tx = \tilde{z}$, $\inner{T^T \ty}{\tx} = 0 \}$. Furthermore, substitution of $T \tx = \tilde{z}$ into~\eqref{eq.mtilde} yields $\tm = \mu \ty$ and~\eqref{eq.tilde} simplifies to, 
	$
	\dot{\tx}
	= 
	-T^T\ty,
	$
	$
	\dot{\ty} 
	= 
	0.
	$	
For~\eqref{eq.tilde}, the largest invariant set $\Omega \DefinedAs \{ (x,y); \nabla f(x) = \nabla f(x^\star)$, $T \tx = \tilde{z}$, $T^T \ty = 0 \}$ $\subseteq \cal D$ is obtained from
	\[
	\inner{T^T \ty}{\tx} 
	\, \equiv \, 
	0
	~ \Rightarrow \;
	\inner{T^T \dot{\ty}}{\tx} \, + \, \inner{T^T \ty}{\dot{\tx}} 
	\, = \,
	- \norm{T^T \ty}^2
	\, \equiv \, 0
	\]
and LaSalle's invariance principle implies that $\Omega$ is GAS.

To complete the proof, we need to show that any $x$ and $y$ that lie in $\Omega$ also satisfy optimality conditions~\eqref{pr.cond} of problem~\eqref{pr.split} with $z = \zms(x,y) = \prox_{\mu g}(Tx + \mu y)$ and thus that $x$ solves problem~\eqref{pr}. For any $(x,y) \in \Omega$, $\nabla f(x) = \nabla f(x^\star)$ and $T^T y = T^T y^\star$. Since $x^\star$ and $y^\star$ are optimal primal-dual points, we have
	\[
	\nabla f(x) \,+\, T^T y
	\; = \;
	\nabla f(x^\star) \,+\, T^T y^\star
	\; = \; 0
	\]
which implies that every $( x, y ) \in \Omega$ satisfies~\eqref{pr.condgrad}. Optimality condition~\eqref{pr.condsubgrad} for $(x^\star, y^\star)$, $Tx^\star = z^\star$, together with $T\tilde{x} = \tilde{z}$, imply that $Tx = z$, i.e., $x$ and $z = \prox_{\mu g}(Tx + \mu y)$ satisfy~\eqref{pr.condfeas}. Finally, the optimality condition of the problem~\eqref{eq.prox} that defines $\prox_{\mu g}(v)$ is
$
	{\partial g(z) + \tfrac{1}{\mu}(z - v) \ni 0.}
$
Letting $v = Tx + \mu y$ from the expression~\eqref{eq.zstar} that characterizes the $\zms$-manifold and noting $Tx = z$ by~\eqref{pr.condfeas} leads to~\eqref{pr.condsubgrad}. Thus, every $(x,y) \in \Omega$ satisfies~\eqref{pr.cond}, implying that the set of primal-dual optimal points is GAS. 
\end{IEEEproof}
	
	\vspace*{-2ex}
\subsection{Global exponential stability for strongly convex $f$}
	\label{sec.IQCs}
	
We express~\ref{eq.ctdyn}, or equivalently~\ref{eq.ctdyn1}, as a linear system $G$ connected in feedback with nonlinearities that correspond to the gradients of $f$ and of the Moreau envelope of $g$; see Fig.~\ref{fig.iqc1}. These  nonlinearities can be conservatively characterized by IQCs. Exponential stability of $G$ connected in feedback with {\em any\/} nonlinearity that satisfies these IQCs implies exponential convergence of~\ref{eq.ctdyn} to the primal-dual optimal solution of~\eqref{pr.split}. In what follows, we adjust the tools of~\cite{lesrecpac16,husei16} to our setup and establish global exponential stability by evaluating the feasibility of an LMI. We assume that the function $f$ is $m_f$-strongly convex with an $L_f$-Lipschitz continuous gradient. Characterizing additional structural restrictions on $f$ and $g$ with IQCs may lead to tighter bounds on the rate of convergence.

As illustrated in Fig.~\ref{fig.iqc1},~\ref{eq.ctdyn1} can be expressed as a linear system $G$ connected via feedback to a nonlinear block $\Delta$,
\[
	\ba{c}
        \dot w
        \; = \;
        A \, w
        \;+\;
        B \, u,
        ~~~
        \xi
        \; = \;
        C \, w,
        ~~~
        u
        \;=\;
        \Delta(\xi)
	\\[0.15cm]
	A
	=
	\tbt{\!\!-m_f I\!\!}{\!\!\!\!}{\!\!\!\!}{\!\!-\mu I\!\!},
	~
	B
	=
	\tbt{\!\!-I\!\!}{\!\!-\frac{1}{\mu}T^T\!\!}{\!\!0\!\!}{\!\!I\!\!},
	~
	C
	=
	\tbt{\!\!I\!\!}{\!\!0\!\!}{\!\!T\!\!}{\!\!\mu I\!\!}
	\ea
\]
where $w \DefinedAs [ \, x^T \; y^T]^T$, $\xi \DefinedAs [ \, \xi_1^T \; \xi_2^T ]^T$, and 
	$
	u \DefinedAs [ \, u_1^T \; u_2^T ]^T.
	$
Nonlinearity $\Delta$ maps the system outputs $\xi_1 = x$ and $\xi_2 = Tx + \mu y$ to the inputs $u_1$ and $u_2$ via $u_1 = \Delta_1(\xi_1) \DefinedAs \nabla f(\xi_1) - m_f \xi_1$ and $u_2 = \Delta_2(\xi_2) \DefinedAs \mu\nabla M_{\mu g}(\xi_2) = \xi_2 - \prox_{\mu g}(\xi_2)$.
 
When the mapping $u_i = \Delta_i (\xi_i)$ is the $L_i$-Lipschitz continuous gradient of a convex function, it satisfies the IQC~\cite[Lemma 6]{lesrecpac16}
\be
	\tbo{\xi_i \,-\, \xi_0}{u_i \,-\, u_0}^T
	\tbt{0}{\hat{L}_i I}{\hat{L}_i I}{-2I}
	\tbo{\xi_i \,-\, \xi_0}{u_i \,-\, u_0}
	\, \geq \,
	0
    \label{eq.iqcf}
\ee
where $\hat{L}_i \geq L_i$, $\xi_0$ is some reference point, and $u_0 = \Delta_i (\xi_0)$. Since $f$ is $m_f$-strongly convex, the mapping $\Delta_1 (\xi_1)$ is the gradient of the convex function $f(\xi_1) - (m_f/2) \norm{\xi_1}^2$. Lipschitz continuity of $\nabla f$ with parameter $L_f$ implies Lipschitz continuity of $\Delta_1(\xi_1)$ with parameter $L_1 \DefinedAs L_f - m_f$; thus, $\Delta_1$ satisfies~\eqref{eq.iqcf} with $\hat{L}_1 \geq L_1$. Similarly, $\Delta_2(\xi_2)$ is the scaled gradient of the convex Moreau envelope and is Lipschitz continuous with parameter $1$; thus, $\Delta_2$ also satisfies~\eqref{eq.iqcf} with $\hat{L}_2 \geq 1$. These two IQCs can be combined into
	\be
	(\eta \, - \, \eta_0)^T
	\Pi 
	\, 
	(\eta \, - \, \eta_0) \, \geq \, 0,
	~~~
	\eta 
	\; \DefinedAs \;
	[ \, \xi^T \; u^T ]^T.
	\label{eq.IQC}
	\ee

\begin{figure}[t!]
\begin{center}
\resizebox{0.65\columnwidth}{!}{\input{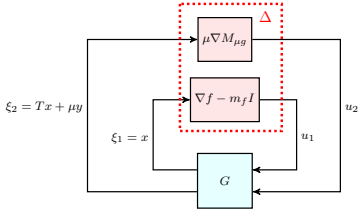}%
\noindent%
%   set a filename for externalization
%\tikzsetnextfilename{clp_Ndelta_config}
%
\begin{tikzpicture}[scale=1, auto, >=stealth']

    \small

   % specify with respect to the plant

     \node[clpblock, minimum height=1.25cm, minimum width=1.25cm] (N)  {$G$};

     \node[deltablock, minimum height=1cm, minimum width=1cm, anchor=south]
                      (delta) at ($(N.north) + (0,0.7cm)$) {$\nabla f - m_fI$};

     \node[deltablock, minimum height=1cm, minimum width=1cm, anchor=south]
                      (delta2) at ($(delta.north) + (0,0.35cm)$) {$\mu\nabla M_{\mu g}$};

     \def\ioOffset{0.25cm}
     \coordinate (Nv) at ($(N.east) + (0,\ioOffset)$);
     \coordinate (Nu) at ($(N.east) - (0,\ioOffset)$);
     \coordinate (Nz) at ($(N.west) + (0,\ioOffset)$);
     \coordinate (Ny) at ($(N.west) - (0,\ioOffset)$);

    % now link the nodes
%    \draw [connector] ($(Nu) + (1.5cm,0)$) -- node [pos=0.1, above] {$d$}
%                                                                node (Numid) [midway] {} (Nu);
%    \draw [connector] (Ny) -- node (Nymid) [midway] {}
%                                         node [pos=0.9, above] {$z$}($(Ny)-(1.5cm,0)$) ;

    \node (Numid) at ($(Nu) + (1.cm,0)$){};
    \node (Numid2) at ($(Numid) + (1.0cm,0)$){};
    \node (Nymid) at ($(Ny)-(1.cm,0)$){};
    \node (Nymid2) at ($(Nymid)-(1.5cm,0)$){};

    \draw [connector] (Nz) -- (Nz -| Nymid) |- (delta);
    \draw [connector]  (delta) -- (delta -| Numid) |- (Nv);

    \draw [connector] (Ny) -- (Ny -| Nymid2) |- (delta2);
    \draw [connector]  (delta2) -- (delta2 -| Numid2) |- (Nu);

    \path (N.north) -- node (midpoint) [midway]{} (delta.south) {};
    \path (N.north) -- node (midpoint2) [midway]{} (delta2.south) {};

    \node[left] at (Nymid |- midpoint) {$\xi_1 = x$};
    \node[left] at (Nymid2 |- midpoint2) {$\xi_2 = Tx + \mu y$};
%    \node[right] at (Numid2 |- midpoint2) {\phantom{$Tx + \mu y$}};
    \node[right] at (Numid |- midpoint) {$u_1$};
    \node[right] at (Numid2 |- midpoint2) {$u_2$\phantom{$= Tx + \mu y$}};

    \draw[red,ultra thick,dotted] ($(delta2.north west)+(-0.4,0.3)$)  rectangle ($(delta.south east)+(0.5,-0.2)$);
    \node [red] at ($(delta2.north east)+(0.3,0.0)$) {\large $\Delta$};

\end{tikzpicture} }
\caption{{Block diagram of gradient flow dynamics~\ref{eq.ctdyn1} where $G$ is a linear system in feedback with nonlinearities that satisfy~\eqref{eq.IQC}.}}
\label{fig.iqc1}
\end{center}
\end{figure}
{For a linear system $G$ connected in feedback with nonlinearities that satisfy IQC~\eqref{eq.IQC},~\cite[Theorem 3]{husei16} establishes $\rho$-exponential convergence, i.e., $\norm{w(t) - w^\star} \leq \tau\mre^{-\rho t}\norm{w(0) - w^\star}$ for some $\tau, \rho > 0$, by verifying the existence of a matrix $P \succ 0$ such that,
\be
	\!\!\!
	\tbt{A_\rho^TP + PA_\rho}{PB}{B^TP}{0}
	\, + \,
	\tbt{C^T}{0}{0}{I} 
	\Pi 
	\tbt{C}{0}{0}{I}
	\, \preceq \,
	0,
	\label{eq.iqclmi}
\ee
where $A_\rho \DefinedAs A + \rho I$. In Theorem~\ref{thm.thmct}, we determine a scalar condition that ensures global exponential stability when $TT^T$ is full rank.
	\begin{theorem} 
	\label{thm.thmct}
Let $f$ be strongly convex with parameter $m_f$, let its gradient be Lipschitz continuous with parameter $L_f$, let $g$ be proper, lower semicontinuous, and convex, and let $TT^T$ be full rank.
Then, if $\mu \geq L_f - m_f$, there is a $\rho > 0$ such that the dynamics~\ref{eq.ctdyn} converge $\rho$-exponentially to the optimal point of~\eqref{pr.split}. 
\end{theorem}}

\begin{IEEEproof}
\begin{subequations}
\label{eq.ctkyp}
Since any function that is Lipschitz continuous with parameter $L$ is also Lipschitz continuous with parameter $\hat{L} \geq L$, we establish the result for $\mu = \hat{L}_1 \DefinedAs L_f - m_f$ and $\hat{L}_2 = 1$. We utilize~\cite[Theorem~3]{husei16} to show $\rho$-exponential convergence by verifying matrix inequality~\eqref{eq.iqclmi} through a series of equivalent expressions~\eqref{eq.ctkyp}. We first apply the KYP Lemma~\cite[Theorem~1]{ran96} to~\eqref{eq.iqclmi} to obtain an equivalent frequency domain characterization
\be
	\tbo{G_\rho (\mrj \omega)}{I}^*
	\Pi
	\tbo{G_\rho (\mrj \omega)}{I}
	\;\preceq\;
	0,
	~~~
	\forall \; \omega \;\in\;\bbR
	\label{eq.kyp}
\ee
where $G_\rho (\mrj \omega) = C(\mrj \omega I - A_\rho)^{-1}B$. Evaluating the left-hand side of~\eqref{eq.kyp} for $L = \mu$ and dividing by $-2$ yields the matrix inequality
\be
	\label{eq.eval}
	\tbt{\!\!\!\dfrac{\mu\hat m + \hat m^2 + \omega^2}{\hat m^2 + \omega^2} \, I\!\!\!}
	{\!\!\!\dfrac{\hat m}{\hat m^2 + \omega^2} \,T^T\!\!\!}
	{\!\!\!*\!\!\!}
	{\!\!\!\dfrac{\hat m/\mu}{\hat m^2 + \omega^2} \,TT^T + \dfrac{\omega^2 - \rho \hat \mu}{\hat\mu^2 + \omega^2} \, I\!\!\!}
	\,\succ\,
	0
\ee
where $\hat m \DefinedAs m_f - \rho > 0$ and $\hat \mu \DefinedAs \mu - \rho > 0$ so that $A_\rho$ is Hurwitz, i.e., the system $G_\rho$ is stable. Since the $(1,1)$ block in~\eqref{eq.eval} is positive definite for all $\omega$, the matrix in~\eqref{eq.eval} is positive definite if and only if the corresponding Schur complement is positive definite,
\be
	\dfrac{ \hat m/\mu}{\mu \hat m + \hat m^2 + \omega^2}
	\, 
	TT^T
	\;+\;
	\dfrac{\omega^2 - \rho \hat \mu}{\hat\mu^2 + \omega^2}
	\,
	I
	\;\succ\;
	0.
	\label{eq.schurLmu}
\ee
We exploit the symmetry of $TT^T$ to diagonalize~\eqref{eq.schurLmu} via a unitary coordinate transformation. This yields $m$ scalar inequalities parametrized by the eigenvalues $\lambda_i$ of $TT^T$. Multiplying the left-hand side of these inequalities by the positive quantity $(\hat \mu^2 + \omega^2)(\mu \hat m + \hat m^2 + \omega^2)$ yields a set of equivalent, quadratic in $\omega^2$, conditions,
\be
	\omega^4
	\,+\,
	( 
	\tfrac{\hat m\lambda_i}{\mu} + \hat m^2 + \mu\hat m - \rho\hat \mu
	)
	\omega^2
	\,+\,
	\hat{m} \hat{\mu}
	(
	\tfrac{\hat{\mu} \lambda_i}{\mu}
	\,-\,
	\rho (\mu + \hat{m})
	) 
	\,>\,
	0.
	\label{eq.quad}
\ee
Condition~\eqref{eq.quad} is satisfied for all $\omega \in \bbR$ if there are no $\omega^2 \geq 0$ for which the left-hand side is nonpositive. When $\rho = 0$, both the constant term and the coefficient of $\omega^2$ are strictly positive, which implies that the roots of~\eqref{eq.quad} as a function of $\omega^2$ are either not real or lie in the domain $\omega^2 < 0$, which cannot occur for $\omega \in \bbR$.  Finally, continuity of~\eqref{eq.quad} with respect to $\rho$ implies the existence a positive $\rho$ that satisfies~\eqref{eq.quad} for all $\omega \in \bbR$.
\end{subequations}
\end{IEEEproof}

\begin{remark}
Each eigenvalue $\lambda_i$ of a full rank matrix $TT^T$ is positive and hence to estimate the exponential convergence rate it suffices to check~\eqref{eq.quad} only for the smallest $\lambda_i$. A sufficient condition for~\eqref{eq.quad} to hold for each $\omega \in \bbR$ is positivity of the constant term and the coefficient multiplying $\omega^2$. For $\rho < \min \, (m_f, \mu)$ these can be, respectively, expressed as the following quadratic \mbox{inequalities in $\rho$,}
	\[
	\ba{rcl}
	\rho^2 
	\, - \,
	\gamma \, \rho
	\, + \,
	\lambda_{\min}
	& \!\!\! > \!\!\! &
	0
	\\[0.05cm]
	2 \rho^2 
	\, - \,
	( \gamma \, + \, \mu \, + \, m_f ) \, \rho
	\, + \,
	\gamma \, m_f 	
	& \!\!\! > \!\!\! &
	0
	\ea
	\] 
where $\gamma \DefinedAs \mu + m_f + \tfrac{\lambda_{\min}}{\mu}$. The solutions to these provide the following estimates of the exponential convergence rate: (i) $\rho < \rho_1$ when $m_f \geq \mu$; and (ii) $\rho  < \min \, (\rho_1, \rho_2)$ when $m_f < \mu$, where	
	\[
	\ba{rcl}
	\rho_1 
	& \!\!\! = \!\!\! &
	\tfrac{1}{2}
	\,
	(
	\gamma
	\, - \,
	\sqrt{\gamma^2 \, - \, 4 \lambda_{\min}}
	)
	\\[0.1cm]
	\rho_2 
	& \!\!\! = \!\!\! &
	\tfrac{1}{4}
	\,
	(
	\gamma \, + \, \mu \, + \, m_f 
	\, - \,
	\sqrt{
	( \gamma \, + \, \mu \, + \, m_f )^2 
	\, - \,
	8 \, \gamma \, m_f 
	}
	).
	\ea
	\]
Our explicit analytical expressions can be used to determine the optimal value of $\mu \geq L_f - m_f$ to maximize the above decay rates.	
	\end{remark}	
\begin{remark}
A similar convergence rate result can be obtained by applying~\cite[Theorem 4]{lesrecpac16} to a discrete-time implementation of the primal-descent dual-ascent dynamics that results from a forward Euler discretization of~\ref{eq.ctdyn}; for details, see~\cite{dinhudhijovCDC18}.
	\end{remark}	
	
	\begin{remark}
To the best of our knowledge, we are the first to establish global exponential stability of the primal-dual gradient flow dynamics for nonsmooth composite optimization problems~\eqref{pr} with a strongly convex $f$. Recent reference~\cite{quli18} proves similar result for a narrower class of problems (strongly convex and smooth objective function with either affine equality or inequality constraints). Both of these can be cast as~\eqref{pr} via introduction of suitable indicator functions and exponential stability follows immediately from our Theorem~\ref{thm.thmct}. This demonstrates power and generality of the proposed approach for nonsmooth composite optimization. While we employ frequency domain IQCs in the proof of Theorem~\ref{thm.thmct}, a time domain Lyapunov-based analysis was used in~\cite{quli18}, which is of \mbox{independent interest.}	
	\end{remark}

	\vspace*{-3ex}
\subsection{Distributed implementation}
	\label{sec.di}
	
Gradient flow dynamics~\ref{eq.ctdyn} are convenient for distributed implementation. If the state vector $x$ corresponds to the concatenated states of individual agents, $x_i$, the sparsity pattern of $T$ and the structure of the gradient map $\nabla f$: $\bbR^n \to \bbR^n$ dictate the communication topology required to form $\nabla\cL_\mu$ in~\ref{eq.ctdyn}. For example, if $f(x) = \sum f_i(x_i)$ is separable over the agents, then $\nabla f_i(x_i)$ can be formed locally. If in addition $T^T$ is an incidence matrix of an undirected network with the graph Laplacian $T^TT$,  each agent {need only} share its state $x_i$ with its neighbors and maintain dual variables $y_i$ that correspond to its edges. {A distributed implementation is also natural when the mapping $\nabla f$: $\bbR^n \to \bbR^n$ is sparse.}

Our approach provides several advantages over existing distributed optimization algorithms. Even for problems~\eqref{pr} with non-differentiable regularizers $g$, a formulation based on the proximal augmented Lagrangian yields gradient flow dynamics~\ref{eq.ctdyn} with a continuous right-hand side. This is in contrast to existing approaches which employ subgradient methods~\cite{nedozdTAC09} or use discontinuous projected dynamics~\cite{feipagAUT10,chemalcorSCL16,cheghacor17,chemallowcor18}. Note that although the augmented Lagrangian $\cL_\mu(x,y;z)$ contains a quadratic term $\tfrac{1}{2\mu}\norm{\cT(x) - z}^2$, it is not {\em jointly\/} strongly convex in $x$ and $z$ and the resulting proximal augmented Lagrangian~\eqref{eq.alprox} is not strictly convex-concave in $x$ and $y$. Furthermore, when $T$ is not diagonal, a distributed proximal gradient cannot be implemented because the proximal operator of $g(Tx)$ may not be separable. Finally, ADMM has been used for distributed implementation in the situations where $f$ is separable and $T$ is an incidence matrix. Relative to such a scheme, our method does not require solving an $x$-minimization subproblem in each iteration and provides a guaranteed rate of convergence.

\begin{remark}
Special instances of our framework have strong connections with the existing methods for distributed optimization on graphs; e.g.,~\cite{nedozdTAC09,waneli11,ghacor14}. The networked optimization problem of minimizing $f(\bar x) = \sum f_i(\bar x)$ over a single variable $\bar x$ can be reformulated as $\sum f_i(x_i) + g(Tx)$ where the components $f_i$ of the objective function are distributed over independent agents $x_i$, $x$ is the aggregate state, $T^T$ is the incidence matrix of a strongly connected and balanced graph, and $g$ is the indicator function associated with the set $Tx = 0$. The $g(Tx)$ term ensures that at feasible points, $x_i = x_j = \bar x$ for all $i$ and $j$. It is easy to show that
	$
    \nabla M_{\mu g}(v)
    =
    (1/\mu) \, v
	$
and that the dynamics~\ref{eq.ctdyn} are given by,
\be
	\ba{rcl}
	\dot x
	&\!\!=\!\!&
	-\, \nabla f(x)
	\;-\;
	(1/\mu) \, L \, x
	\;-\;
	\tilde y
	\\[0.05cm]
	\dot{\tilde y}
	&\!\!=\!\!&
	\beta \, L \, x
	\ea
	\label{eq.gf-network}
\ee
where $\beta > 0$, $L \DefinedAs T^TT$ is the graph Laplacian of a connected undirected network, and $\tilde{y} \DefinedAs T^T y$ belongs to the orthogonal complement of the vector of all ones. The only difference relative to~\cite[Eq.~(20)]{waneli11} and~\cite[Eq.~(11)]{ghacor14} is that $-\tilde{y}$ appears instead of $- L \tilde y$ in equation~\eqref{eq.gf-network} for the dynamics of the primal variable~$x$.
	\end{remark}

	\begin{remark}
Forward Euler discretization of~\eqref{eq.gf-network} is given by
	\be
	\ba{rcl}
	x^{k+1}
	&\!\!=\!\!&
	\left(
	I \, - \, (\alpha/\mu) L
	\right)
	x^k
	\, -\, 
	\alpha \, \nabla f(x^k) 
	\, - \, 
	\alpha
	\,
	\tilde{y}^k
	\\[0.1cm]
	\tilde y^{k+1}
	&\!\!=\!\!&
	\tilde y^k
	\;+\;
	\alpha
	\,
	\beta
	\, 
	L 
	\,
	x^k
	\ea
	\label{eq.gf-network1}
\ee
where $\alpha$ is the step-size, and the EXTRA algorithm~\cite[Equation (2.13)]{shilinyin15}, which has received significant recent attention,
	\be
	x^{k+1}
	\; = \;
	Wx^k
	\; - \;
	\alpha \, \nabla f(x^k)
	\; + \;
	\dfrac{1}{2}
	\;
	{\ds\sum^{k - 1}_{i \, = \, 0} \, (W \, - \, I) \, x^i}
	\label{eq.extra}
	\ee
can be clearly recovered from~\eqref{eq.gf-network1} by setting $\beta = 1/(2 \alpha \mu)$ and taking $W = I - (\alpha/\mu) L$ in~\eqref{eq.extra}. 
\end{remark}
	
	\vspace*{-2ex}
\section{Examples}
\label{sec.ex}

We solve the problems of edge addition in directed consensus networks and optimal placement to illustrate the effectiveness of the proximal augmented Lagrangian method.

	\vspace*{-2ex}
\subsection{Edge addition in directed consensus networks}

A consensus network with $N$ nodes converges to the average of the initial node values $\bar \psi = (1/N) \sum_i \psi_i(0)$ if and only if it is strongly connected and balanced~\cite{mesege10}. Unlike for undirected networks~\cite{wujovSCL17,mogjovTCNS17}, the problem of edge addition in directed consensus networks is not known to be convex. The steady-state variance of the deviations from average is given by the square of the $\cH_2$ norm of,
\be
    \ba{rl}
	\dot \psi
	\; = \; 
	-( L_p + L_x ) \, \psi
	\; + \;
	d,
        &
    \xi
    \; = \; 
    \tbo{Q^{1/2}}{- R^{1/2} L_x} \psi
    \ea
    \non
\ee
where $d$ is a disturbance, $L_p$ is a weighted directed graph Laplacian of a plant network, $Q \DefinedAs I - (1/N)\one\one^T$ penalizes the deviation from average, and $R \succ 0$ is the control weight. The objective is to optimize the $\cH_2$ norm (from $d$ to $\xi$) by adding a few additional edges, specified by the graph Laplacian $L_x$ of a controller network.

To ensure convergence of $\psi$ to the average of the initial node values, we require that the closed-loop graph Laplacian, $L = L_p + L_x$, is balanced. This condition amounts to the linear constraint, $\one^TL = 0$. We express the directed graph Laplacian of the controller \mbox{network as,}
	$
    L_x
    =
    \sum_{i \, \neq \, j}
    L_{ij}
    z_{ij}
    \AsDefined
    \sum_l
    L_l
    z_{l}
$
where $z_{ij} \geq 0$ is the added edge weight that connects node $j$ to node $i$, $L_{ij} \DefinedAs \mre_i\mre_i^T - \mre_i\mre_j^T$, $\mre_i$ is the $i$th basis vector in $\bbR^n$, and the integer $l$ indexes the edges such that $z_l = z_{ij}$ and $L_l = L_{ij}$. For simplicity, we assume that the plant network $L_p$ is balanced and connected. Thus, enforcing that $L$ is balanced amounts to enforcing {the linear constraint}
$
	{
    \one^TL_x
    =
    \one^T
    (
    \sum
    L_l
    z_l
    )
    \AsDefined
    (Ez)^T
    =
    0}
$
on $z$, where $E$ is the incidence matrix~\cite{mesege10} of the edges that may be added. Any vector of edge weights $z$ that {satisfies this constraint} can be written as 
$z = Tx$ where the columns of $T$ span the nullspace of the matrix $E$ and provide a basis for the space of balanced graphs, i.e., the cycle space~\cite{mesege10}. Each feasible controller {Laplacian} can thus be written as,
	\begin{subequations}
\be
    L_x
    \, = \,
    \sum_{l}
    L_{l}
    \,
    [ Tx ]_{l}
    \, = \,
    \sum_{l}
    L_{l}
    \bigg[
    \ds\sum_k
    \,
    (T\mre_k)
    \,
    x_k
    \bigg]_{l}
   \, \AsDefined \,
    \ds\sum_{k}
    \hat L_k
    \,
    x_k
    \label{eq.Lx}
	\ee
where the matrices $\hat L_k$ are given by
$
    \hat L_k
    =
    \sum_{l}
    L_{l}
    \,
    [ T\mre_k ]_{l}.
    \label{eq.hatE}
$

Since the mode {corresponding to $\one$} is marginally stable, unobservable, and uncontrollable, we introduce a change of coordinates to the deviations from average $\phi = V^T\psi$ where $V^T\one = 0$ and discard the average mode $\bar \psi = \one^T\psi$. The energy of the deviations from average is given by the the $\Ht$ norm squared of the reduced system,
\be
    f(x)
    \, = \,
    \inner{V^T(Q \, + \, L_x^TRL_x) \, V}{X},
    ~
    \hat A  \, X + X \hat A^T + \hat B \hat B^T 
    \,=\,
    0
    \label{eq.ht}
\ee
\end{subequations}
where $X$ is the controllability gramian of the reduced system with $\hat A \DefinedAs - V^T(L_p + L_x)V$ and $\hat B \DefinedAs V^T$.

To balance the closed-loop $\cH_2$ performance with the number of added edges, we introduce a regularized optimization problem
\be
    x_\gamma
    \;=\;
    \argmin\limits_x
    ~ \;
    f(x)
    \;+\;
    \gamma \, \one^TTx
    \;+\;
    I_+(Tx).
    \label{pr.dir}
\ee
Here, the regularization parameter $\gamma > 0$ specifies the emphasis on sparsity relative to the closed-loop performance $f$, and $I_+$ is the indicator function associated with the nonnegative orthant $\bbR^m_+$. When the desired level of sparsity for the vector of the added edge weights $z_\gamma = Tx_\gamma$ has been attained, optimal weights for the identified set of edges are obtained by solving,
\be
    \ba{rl}
    \minimize\limits_x
    &
    f(x)
    \;+\;
    I_{\cal Z_\gamma}(Tx)
    \;+\;
    I_+(Tx)
    \ea
    \label{pr.dirpol}
\ee
where $\cal Z_\gamma$ is the set of vectors with the same sparsity pattern as $z_\gamma$ and $I_{\cal Z_\gamma}$ is the indicator function associated with this set.

\subsubsection{Implementation}

We next provide implementation details for solving~\eqref{pr.dir} and~\eqref{pr.dirpol}. The proof of next lemma is omitted for brevity.

\begin{lemma}
Let a graph Laplacian of a directed plant network $L_p$ be balanced and connected and let $\hat A$, $\hat B$, $L_x$, and $V$ be as defined in~\eqref{eq.Lx}--\eqref{eq.ht}. The gradient of $f(x)$ defined in~\eqref{eq.ht} is given by,
\[
    \nabla f(x)
    \;=\;
    2
    \vect
    \left( \inner{
    (
    R \, L_x V
    \, - \,
    VP
    )
    \, XV^T}{\hat L_k}
    \right)
\]
\mbox{where $X$ and $P$ are the} controllability and observability gramians determined by~\eqref{eq.ht} and $\hat A^TP + P\hat A + V^T(Q + L_x^TRL_x)  V = 0$.
\end{lemma}

The proximal operator associated with the regularization function \mbox{{$g_s (z) \DefinedAs \gamma \one^T z + I_+(z)$}} in~\eqref{pr.dir} is 
$
    \prox_{\mu g_s}(v_i)
    =
    \max
    \{
    0,
    v_i -
    \gamma \mu
    \},
$
the Moreau envelope is given by
$
    M_{\mu g_s}(v)
    =
	\sum_i
    \{  
    v_i^2/(2\mu),
    \,
    v_i \leq \gamma \mu; 
    \,
    \gamma \, (v_i - \gamma \mu/2 ),
    \,
    v_i > \gamma \mu
    \},
$
and
$
	\nabla M_{\mu g_s}(v)
	=
	\max \, \{v/\mu, \, \gamma\}.
$
The proximal operator of the regularization function in~\eqref{pr.dirpol}, $g_p (z) \DefinedAs I_{\cal Z_\gamma} (z) + I_+(z)$, is a projection onto the intersection of the set $\cal Z_\gamma$ and the nonnegative orthant,
$
	\prox_{\mu g_p}(v)
	=
	\cP_{\cE}(v),
	$
the Moreau envelope is the distance to $\cE \DefinedAs {\cal Z_\gamma} \cap \bbR^m_+$,
$
	M_{\mu g_p}(v)
	=
	\tfrac{1}{2\mu}
	\norm{v \,-\, \cP_{\cE}(v)}^2
$
and $\nabla M_{\mu g_p} (v)$ is determined by a vector pointing from $\cE$ to $v$,
$
	\nabla M_{\mu g_p}(v)
	=
	\tfrac{1}{\mu} 
	(
	v - \cP_{\cE}(v)).
$

\subsubsection{Computational experiments}

We solve~\eqref{pr.dir} and~\eqref{pr.dirpol} using Algorithm~\ref{alg.mm}, where L-BFGS is employed in the $x$-minimization subproblem~\ref{eq.MMa1}. For the plant network shown in Fig.~\ref{fig.ex2}, Fig.~\ref{fig.perfvsp} illustrates the tradeoff between the number of added edges and the closed-loop $\Ht$ norm. The added edges are identified by computing the $\gamma$-parameterized homotopy path for problem~\eqref{pr.dir}, and the optimal edge weights are obtained by solving~\eqref{pr.dirpol}. The red dashed lines in Fig.~\ref{fig.ex2} show the optimal set of $2$ added edges. These are obtained for $\gamma = 3.5$ and they yield $23.91\%$ performance loss relative to the setup in which all edges in the controller graph are used. We note that the same set of edges is obtained by conducting an exhaustive search. This suggests that the proposed convex regularizers may offer a good proxy for solving difficult combinatorial optimization problems.

We also consider simple directed cycle graphs with $N = 5$ to $50$ nodes and $m = N^2 - N$ potential added edges. We solve~\eqref{pr.dir} for $\gamma = 0.01, 0.1, 0.2$, and $R = I$ using the proximal augmented Lagrangian MM algorithm (PAL), ADMM, and ADMM with an adaptive heuristic for updating $\mu$~\cite{boyparchupeleck11} (ADMM $\mu$). The $x$-update in each algorithm is obtained using L-BFGS. Since $g_s(Tx)$ and $g_p(Tx)$ are not separable in $x$, proximal gradient cannot be used here.

Figure~\ref{fig.scale} shows the time required to solve problem~\eqref{pr.dir} in terms of the total number of potential added edges; Fig.~\ref{fig.outer} demonstrates that PAL requires fewer outer iterations; and Fig.~\ref{fig.periter} illustrates that the average computation time per outer iteration is roughly equivalent for all three methods. Even with an adaptive update of $\mu$, ADMM requires more outer iterations which increases overall solve time relative to the proximal augmented Lagrangian method. Thus, compared to ADMM, PAL provides computational advantage by reducing the number of outer iterations \mbox{(indexed by $k$ in Algorithm~\ref{alg.mm} and in~\eqref{eq.ADMM}).}

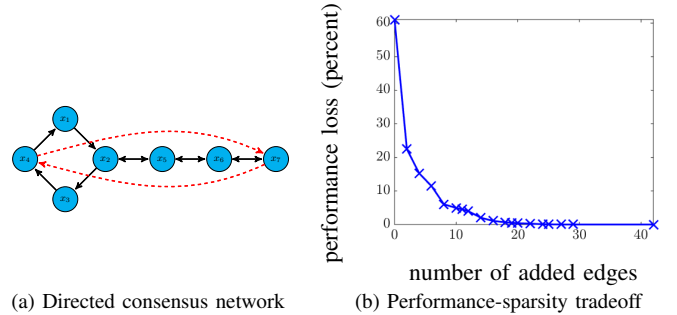
\begin{figure}[b!]
\captionsetup[subfloat]{farskip=2pt,captionskip=1pt}
    \begin{tabular}{cc}
 \subfloat[Directed consensus network\label{fig.ex2}]{
 \begin{tabular}{c}
 \!\!\!\!\!
 	\vspace{1cm}
	\!\!\!\!\!
 	\\
	\!\!\!\!\!
        \resizebox{.42\columnwidth}{!}{\begin{tikzpicture}[>=stealth',shorten >=1pt, node distance=1.8cm, on grid, initial/.style={}]

  \node[state,fill=ProcessBlue]          (1)                        {$x_1$};
  \node[state,fill=ProcessBlue]          (2) [below right =of 1]    {$x_2$};
  \node[state,fill=ProcessBlue]          (4) [below left =of 1]    {$x_4$};
  \node[state,fill=ProcessBlue]          (3) [below left =of 2]    {$x_3$};
  \node[state,fill=ProcessBlue]          (5) [right =of 2]    {$x_5$};
  \node[state,fill=ProcessBlue]          (6) [right =of 5]    {$x_6$};
  \node[state,fill=ProcessBlue]          (7) [right =of 6]    {$x_7$};
%  \node[state,scale=0.7]          (5) [below right =of 4]    {$5$};
%  \node[state,scale=0.7]          (6) [above left =of 4]    {$6$};

\tikzset{every node/.style={fill=white}}

%\tikzset{mystyle/.style={->,double=black}}
%\path
%% (1)     edge [mystyle]    (2)
%      (1)     edge [mystyle]    (3)
%%      (5)     edge [mystyle]    (4)
%%      (3)     edge [mystyle]    (5)
%%      (5)     edge [mystyle]    (6)
%      (3)     edge [mystyle]    (4)
%      (4)     edge [mystyle]    (1);

%\tikzset{mystyle/.style={<->,double=black}}
%     \path (4)     edge [mystyle]    node   {$1$} (1);
%

% curvy arrows
\tikzset{mystyle/.style={->,relative=false,in=135,out=-45,double=black}}
\path (1)     edge [mystyle]   (2);

\tikzset{mystyle/.style={->,relative=false,in=45,out=-135,double=black}}
\path (2)     edge [mystyle]   (3);

\tikzset{mystyle/.style={->,relative=false,in=-45,out=135,double=black}}
\path (3)     edge [mystyle]   (4);

\tikzset{mystyle/.style={->,relative=false,in=-135,out=45,double=black}}
\path (4)     edge [mystyle]   (1);

\tikzset{mystyle/.style={<->,relative=false,in=180,out=0,double=black}}
\path (2)     edge [mystyle]   (5);

\tikzset{mystyle/.style={<->,relative=false,in=180,out=0,double=black}}
\path (5)     edge [mystyle]   (6);

\tikzset{mystyle/.style={<->,relative=false,in=180,out=0,double=black}}
\path (6)     edge [mystyle]   (7);

%\tikzset{mystyle/.style={->,relative=false,in=0,out=-155,dashed,draw=red,double=red}}
%\path (7)     edge [mystyle]   (3);
%
%\tikzset{mystyle/.style={->,relative=false,in=-90,out=90,dashed,draw=red,double=red}}
%\path (3)     edge [mystyle]   (1);
%
%\tikzset{mystyle/.style={->,relative=false,in=155,out=0,dashed,draw=red,double=red}}
%\path (1)     edge [mystyle]   (7);

\tikzset{mystyle/.style={->,relative=false,in=-15,out=-155,dashed,draw=red,double=red}}
\path (7)     edge [mystyle]   (4);

\tikzset{mystyle/.style={->,relative=false,in=155,out=15,dashed,draw=red,double=red}}
\path (4)     edge [mystyle]   (7);

\end{tikzpicture} }\!\!\!\!\!
	\\
	\!\!\!\!\!
	\vspace{.52cm}\!\!\!\!\!
\end{tabular}
    }
    \!\!\!\!\!&\!\!\!\!\!
  \subfloat[Performance-sparsity tradeoff\label{fig.perfvsp}]{
\begin{tabular}{rc}
   \rotatebox{90}{\!\!\!performance loss (percent)}
    \!\!\!\!\!
    &
    \!\!\!\!\!
    {\includegraphics[width=0.5\columnwidth]{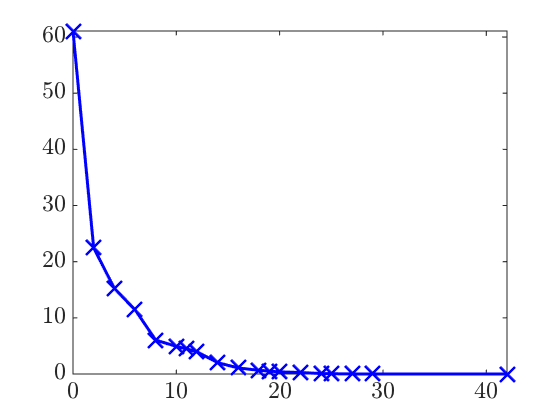}}
    \\
    &
    number of added edges
    \end{tabular}}
    \end{tabular}
\caption{(a) A balanced plant graph with $7$ nodes and $10$ directed edges (solid black lines). A sparse set of $2$ added edges (dashed red lines) is identified by solving~\eqref{pr.dir} with $\gamma = 3.5$ and $R = I$. (b) Tradeoff between performance and sparsity resulting from the solution to~\eqref{pr.dir}-\eqref{pr.dirpol} for the network shown in Fig.~\ref{fig.ex2}. Performance loss is measured relative to the optimal centralized controller (i.e., all edges are used).}
\end{figure}

\vspace*{-3ex}
\subsection{Optimal placement problem}
\vspace{-1ex}

To illustrate the utility of our primal-descent dual-ascent approach, we consider an example in which mobile agents aim to minimize their Euclidean distances relative to a set of targets $\{ b_i \}$ while staying within a desired distance from their neighbors in a network with the incidence matrix $T$,
	\be
	\ba{rl}
	\minimize\limits_x
	&
	\sum_i \, (x_i \,-\, b_i)^2
	\;+\;
	I_{[-1,1]}(Tx).
	\ea
	\label{pr.pddaex}
\ee
Here, $Tx$ is a vector of inter-agent distances which must be kept within an interval $[-1, 1]$. Applying primal-descent dual-ascent update rules to~\eqref{pr.pddaex} achieves path planning for first-order agents $\dot x = u$ with $u = -\nabla_x\cL_\mu(x;y)$. The proximal operator is projection onto a box, $\prox_{\mu I_{[-1,1]}}(z) = \max(\min(z,1),-1)$, the Moreau envelope is the distance squared to that set, $M_{\mu I_{[-1,1]}}(z) = \tfrac{1}{2\mu}\sum\cS_1^2(z_i)$, and $\nabla M_{\mu I_{[-1,1]}} (z) = \tfrac{1}{\mu}\cS_1(z)$. To update its state, each agent $x_i$ needs information from its neighbors in a network with a Laplacian $T^TT$.

Methods based on the subderivative are not applicable because the indicator function is not subdifferentiable. Proximal methods are hindered because the proximal operator of $I_{[-1,1]}(Tx)$ is difficult to compute due to $T$. Since $f(x) = \sum (x_i - b_i)^2$ is separable, a distributed ADMM implementation can be applied; however, it may require large discrete jumps in agent positions, which could be unsuitable for vehicles. Moreover, when $f$ is not separable a distributed implementation of the $x$-minimization step~\eqref{eq.ADMMx} in ADMM would not be possible.

Figure~\ref{fig.pddaex} shows an implementation for a problem with $5$ agents whose set of targets changes position at time $5$. The primal-descent dual-ascent dynamics~\ref{eq.ctdyn} are simulated in {\sc Matlab} using ode$45$.

	\begin{figure*}
  \centering
   \subfloat[Total solve time ($s$) \label{fig.scale}]{
        \begin{tabular}{c}
        {\includegraphics[width=0.2\textwidth]{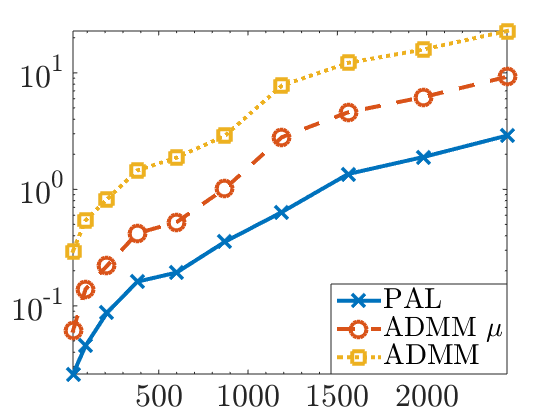}}
        \\
        $m$
        \end{tabular}
    }
    \subfloat[Number of outer iterations \label{fig.outer}]{
        \begin{tabular}{c}
        {\includegraphics[width=0.2\textwidth]{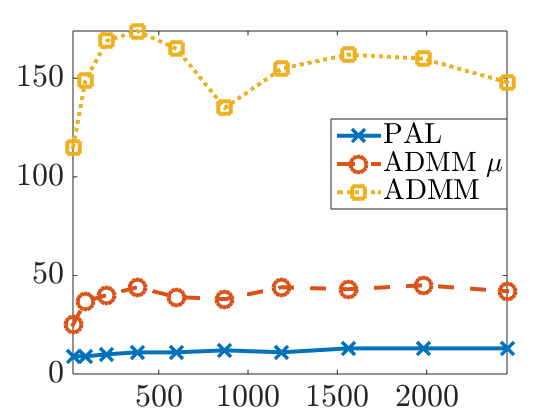}}
        \\
        $m$
        \end{tabular}
    }
    \subfloat[Solve time ($s$) per outer iteration\label{fig.periter}]{%
        \begin{tabular}{c}
        {\includegraphics[width=0.2\textwidth]{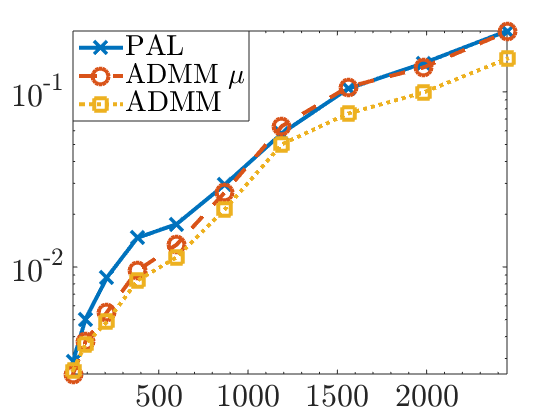}}
        \\
        $m$
        \end{tabular}
    }
    \caption{(a) Total time; (b) number of outer iterations; and (c) average time per outer iteration required to solve~\eqref{pr.dir} with $\gamma = 0.01, 0.1, 0.2$ for a cycle graph with $N = 5$ to $50$ nodes as a function of $m = 20$ to $2450$ potential added edges using PAL (\tc{blue}{\bf --$\bf \times$--}), ADMM (\tcr{\bf - -$\circ$- -}), and ADMM with the adaptive $\mu$-update heuristic~\cite{boyparchupeleck11} (\tc{newyellow}{\bf $\cdots\square\cdots$}). PAL requires fewer outer iterations and thus a smaller total solve time.}
    \label{fig.mor}
  \end{figure*}

	\vspace*{-2ex}
\section{Concluding remarks}
\label{sec.remarks}
	
For a class of nonsmooth composite optimization problems that arise in structured optimal control, we have introduced continuously differentiable proximal augmented Lagrangian function. This function is obtained by collapsing the associated augmented Lagrangian onto the manifold resulting from explicit minimization over the variable in the nonsmooth part of the objective function. Our approach facilitates development of customized algorithms based on the method of multipliers and the primal-descent dual-ascent method.

MM based on the proximal augmented Lagrangian is applicable to a broader class of problems than proximal gradient methods, and it has more robust convergence guarantees, more rigorous parameter update rules, and better practical performance than ADMM. The primal-descent dual-ascent gradient dynamics we propose are suitable for distributed implementation and have a continuous right-hand side. When the differentiable component of the objective function is (strongly) convex, we establish global (exponential) asymptotic stability. Finally, we illustrate the efficacy of our algorithms using the edge addition and optimal placement problems. Future work will focus on developing second-order updates for the primal and dual variables and on providing an extension to nonconvex regularizers. 

\begin{figure}[h!]
\centering
\begin{tabular}{rc}
    \rotatebox{90}{\quad  position}
    \!\!\!\!\!
    &
    \!\!\!\!\!
    {\includegraphics[width=0.75\columnwidth]{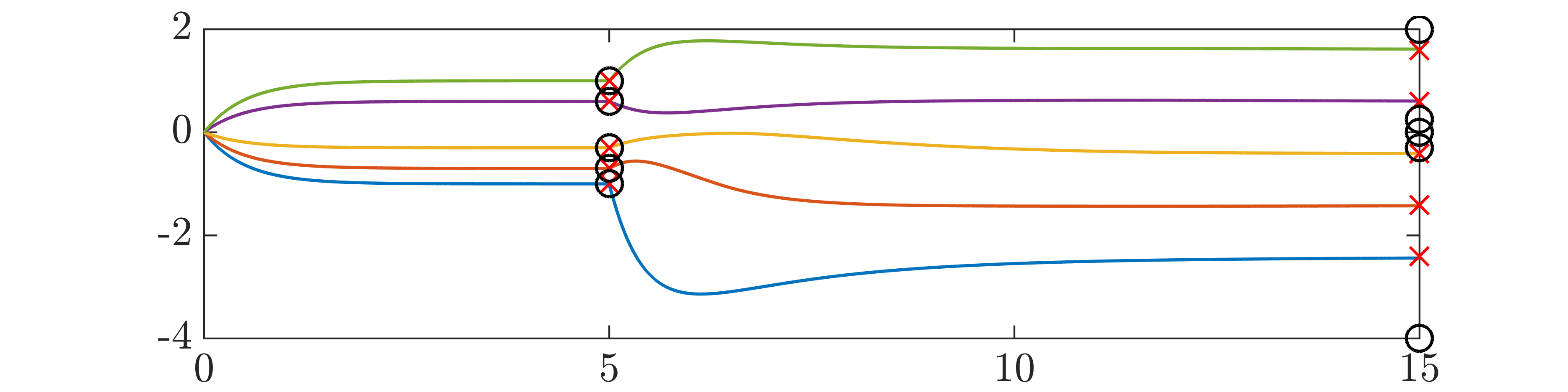}}
    \\
    &
    \vspace{-1ex}
    time
    \end{tabular}
\caption{Set of $5$ distributed agents tracking targets ({\bf black $\circ$}) whose optimal positions are determined by the solution to~\eqref{pr.pddaex} ({\bf \tcr{red $\times$}}).}
\label{fig.pddaex}
\end{figure}

\end{document}